\documentclass[letterpaper,reqno,11pt,oneside,final]{amsart}
\pdfoutput=1
\usepackage{amsmath,amsthm,amsfonts,amssymb}
\usepackage[extension=pdf]{hyperref}
\usepackage{enumitem,comment}
\usepackage[T1]{fontenc}
\usepackage{times}
\usepackage{mathtools}
\usepackage[dvipsnames]{xcolor}
\usepackage{tikz}
\usepackage{mathrsfs}
\usepackage{cases}
\usepackage{bm}
\usepackage[multiple]{footmisc}
\usepackage{microtype}
\usepackage[totalheight=9.58in,totalwidth=6.15in,centering]{geometry}
\usepackage{graphics,graphicx}
\usepackage{tocvsec2}
\usepackage{xr}
\usepackage[style=alphabetic,sorting=nyt,natbib=true,maxnames=99,isbn=false,doi=false,url=false,firstinits=true,hyperref=auto,arxiv=abs,backend=bibtex]{biblatex}
\AtEveryBibitem{%
  \clearlist{language}}
\DeclareFieldFormat[article,inbook,incollection,inproceedings,patent,thesis,unpublished]{title}{#1\isdot}
\renewbibmacro{in:}{%
  \ifentrytype{article}{}{%
    \printtext{\bibstring{in}\intitlepunct}}} 
\defbibheading{apa}[\refname]{\section*{#1}}
\DeclareFieldFormat{sentencecase}{\MakeSentenceCase{#1}} \renewbibmacro*{title}{%
  \ifthenelse{\iffieldundef{title}\AND\iffieldundef{subtitle}}{}
  {\ifthenelse{\ifentrytype{article}\OR\ifentrytype{inbook}%
      \OR\ifentrytype{incollection}\OR\ifentrytype{inproceedings}%
      \OR\ifentrytype{inreference}} {\printtext[title]{%
        \printfield[sentencecase]{title}%
        \setunit{\subtitlepunct}%
        \printfield[sentencecase]{subtitle}}}%
    {\printtext[title]{%
        \printfield[titlecase]{title}%
        \setunit{\subtitlepunct}%
        \printfield[titlecase]{subtitle}}}%
    \newunit}%
  \printfield{titleaddon}}

\definecolor{darkblue}{rgb}{0.13,0.13,0.39}%
\hypersetup{colorlinks=true,urlcolor=darkblue,citecolor=darkblue,linkcolor=darkblue,%
  pdftitle=Kolmogorov equation for TASEP,pdfauthor={M. Nica, J. Quastel, D. Remenik}}

\mathtoolsset{showmanualtags,showonlyrefs}

\newcommand{\fpeqref}[1]{\cite[Eqn. (\ref*{fp-#1})]{fixedpt}}
\newcommand{\fprefst}[1]{\ref*{fp-#1}}

\newtheorem{thm}{Theorem}[section] 

\newtheorem{theorem}{Theorem}[section]

\theoremstyle{definition} \newtheorem{rem}[thm]{Remark} \newtheorem*{rem*}{Remark}

 \newcounter{assum}

\newcommand{\xx}{X}

\newcommand{\R}{R}

\newcommand{\I}{{\rm i}}
\newcommand{\pp}{\mathbb{P}}

\newcommand{\ee}{\mathbb{E}}
\newcommand{\rr}{\mathbb{R}}
\newcommand{\nn}{\mathbb{N}}
\newcommand{\zz}{\mathbb{Z}}

\newcommand{\cL}{\mathcal{L}}

\newcommand{\p}{\partial}
\renewcommand{\d}{\mathrm{d}}
\newcommand{\uno}[1]{\mathbf{1}_{#1}}

\newcommand{\vs}{\vspace{6pt}}
\newcommand{\wt}{\widetilde}

\newcommand{\K}{K_{\Ai}}
\newcommand{\m}{\overline{m}}

\newcommand{\qand}{\quad\text{and}\quad}

\newcommand{\ts}{\hspace{0.1em}}
\newcommand{\tts}{\hspace{0.05em}}
\newcommand{\tsm}{\hspace{-0.1em}}
\newcommand{\ttsm}{\hspace{-0.05em}}

\newcommand{\inv}[1]{\frac{1}{#1}}

\newcommand{\itwopii}[1]{\frac{1}{(2\pi\I)^{#1}}}

\newcommand{\TASEP}{\uptext{TASEP}}

\makeatletter
\newcommand\RedeclareMathOperator{%
  \@ifstar{\def\rmo@s{m}\rmo@redeclare}{\def\rmo@s{o}\rmo@redeclare}%
}
\newcommand\rmo@redeclare[2]{%
  \begingroup \escapechar\m@ne\xdef\@gtempa{{\string#1}}\endgroup
  \expandafter\@ifundefined\@gtempa
     {\@latex@error{\noexpand#1undefined}\@ehc}%
     \relax
  \expandafter\rmo@declmathop\rmo@s{#1}{#2}}
\newcommand\rmo@declmathop[3]{%
  \DeclareRobustCommand{#2}{\qopname\newmcodes@#1{#3}}%
}
\@onlypreamble\RedeclareMathOperator
\makeatother
\newcommand{\uptext}[1]{\text{\upshape{#1}}}
\DeclareMathOperator{\epi}{\uptext{epi}}

\DeclareMathOperator{\Ai}{\uptext{Ai}}
\DeclareMathOperator{\tr}{\mathop{\uptext{tr}}}
\RedeclareMathOperator{\det}{\mathop{\uptext{det}}}
\RedeclareMathOperator{\ker}{\mathop{\uptext{ker}}}
\RedeclareMathOperator{\exp}{\mathop{\uptext{exp}}}
\RedeclareMathOperator{\log}{\mathop{\uptext{log}}}
\RedeclareMathOperator*{\lim}{\mathop{\uptext{lim}}}
\RedeclareMathOperator*{\sup}{\mathop{\uptext{sup}}}
\RedeclareMathOperator*{\limsup}{\mathop{\uptext{lim\hspace{1pt}sup}}}
\RedeclareMathOperator*{\max}{\mathop{\uptext{max}}}
\RedeclareMathOperator*{\inf}{\mathop{\uptext{inf}}}
\RedeclareMathOperator*{\min}{\mathop{\uptext{min}}}

\newcommand{\SM}{\mathcal{S}}
\newcommand{\SN}{\bar{\mathcal{S}}}
\renewcommand{\P}{\chi}

\def\dash---{\kern.16667em---\penalty\exhyphenpenalty\hskip.16667em\relax}

\numberwithin{equation}{section}

\let\oldmarginpar\marginpar
\renewcommand\marginpar[1]{\-\oldmarginpar[\raggedleft\footnotesize #1]%
  {\raggedright{\small\textsf{#1}}}}

\allowdisplaybreaks[2]

\begin{document}

\title{Solution of the Kolmogorov equation for TASEP}

\date{\today}

\author{Mihai Nica} \address[M.~Nica]{
  Department of Mathematics\\
  University of Toronto\\
  40 St. George Street\\
  Toronto, Ontario\\
  Canada M5S 2E4} \email{mnica@math.toronto.edu}

\author{Jeremy Quastel} \address[J.~Quastel]{
  Department of Mathematics\\
  University of Toronto\\
  40 St. George Street\\
  Toronto, Ontario\\
  Canada M5S 2E4} \email{quastel@math.toronto.edu}

\author{Daniel Remenik} \address[D.~Remenik]{
  Departamento de Ingenier\'ia Matem\'atica and Centro de Modelamiento Matem\'atico (UMI-CNRS 2807)\\
  Universidad de Chile\\
  Av. Beauchef 851, Torre Norte, Piso 5\\
  Santiago\\
  Chile} \email{dremenik@dim.uchile.cl}

\begin{abstract}  
We provide a direct and elementary proof that the formula obtained in \cite{fixedpt} for the TASEP transition probabilities for general (one-sided) initial data solves the Kolmogorov backward equation.
The same method yields the solution for the related PushASEP particle system.
\end{abstract}

\maketitle

\section{Introduction}

The \emph{totally asymmetric simple exclusion process} (TASEP) consists of  particles on the lattice $\zz$ performing totally asymmetric nearest neighbour random walks with exclusion: Each particle independently attempts jumps to the neighbouring site to the right at rate $1$, the jump being allowed only if that site is unoccupied.
We will always consider initial conditions in which there is a rightmost particle; this always remains so, and the positions of the particles can be denoted $X_t(1)>X_t(2)>\dotsm$.
The dynamics of the first $N$ particles $X_t(1)>X_t(2)>\dotsm>X_t(N)$ is independent of the rest, so the infinite system clearly makes sense.
If we let
\begin{equation}
\xx^{-1}_t(u) = \min \{k \in \zz : \xx_t(k) \leq u\}
\end{equation}
denote the label of the rightmost particle which sits to the left of, or at, $u$ at time $t$, then the \emph{TASEP  height function} associated to $X_t$ can be defined, for $z\in\zz$, by 
\[h_t(z) = -2\tsm\left(\xx_t^{-1}(z-1) - \xx_0^{-1}(-1) \right) - z,\]
which fixes $h_0(0)=0$.
The height function itself is a simple random walk path $h_t(z+1) = h_t(z) +\hat{\eta}_t(z)$ with $\hat{\eta}_t(z)=1$ if there is a particle at $z$ at time $t$ and $-1$ if there is no particle at $z$ at time $t$.
The dynamics of the height function  is that local max's become local min's at rate $1$; i.e. if $h_t(z) = h_t(z\pm 1) +1$ then $h_t(z)\mapsto h_t(z)-2$ at rate $1$, the rest of the height function remaining unchanged.  The rate of decrease is
\begin{equation}
-2\uno{\wedge} = \tfrac12\left[ ( \nabla^-h)( \nabla^+h) - 1 + \Delta h\right]
\end{equation}
where $\uno{\wedge}$ is the indicator function of a local max.  Hence the TASEP height function can be seen as a 
simple discretization of the Kardar-Parisi-Zhang (KPZ) equation,
\begin{equation}
\partial_t h = \tfrac12 (\partial_x h)^2 +\tfrac12 \partial_x^2 h  + \xi.
\end{equation}
Because of its amenability to computations, TASEP has become the most popular model in the KPZ universality class. 

$N$-particle TASEP was solved by \citet{MR1468391} using the coordinate Bethe ansatz.   
The transition probabilities are given by a determinant,
\begin{equation}\label{eqGreen}
\pp (X_t( 1)= x_1,\ldots,X_t(N)=x_N)=\det(G_{i-j}(t,x_{N+1-i}-X_0(N+1-j)))_{1\leq i,j\leq N}
\end{equation}
with
\begin{equation}\label{eqFn}
G_{n}(t,x)=\frac{(-1)^n}{2\pi \I} \oint_{\Gamma_{0,1}}\d w\,\frac{(1-w)^{-n}}{w^{x-n+1}}e^{t(w-1)},
\end{equation}
where $\Gamma_{0,1}$ is any positively oriented simple loop which includes $w=0$ and $w=1$.
A direct proof of this formula is not difficult and can be obtained in a couple of pages (see, for example \cite{mq-notes}).

On the other hand, one is generally less interested at a later time in the exact positions of the particles, or the exact height function, but in the joint distribution of the height function at a finite number of points $x_1,\ldots, x_m$, where $m$ is fixed, and $N$ is large, or infinite. 
Or, what amounts to the same thing, one would like to compute, for some sequence $n$ of $m$ indices $n_1 < \dotsc < n_m$ and any vector $a = (a_1, \ldots, a_m)\in\zz^m$, the joint probability
\begin{equation}\label{eq:extKernelProb}
 F_t(X_0;a,n)=\pp_{X_0}\!\left(X_t(n_j)>a_j,~j=1,\dotsc,m\right),
\end{equation}
where the subscript in $\pp_{X_0}$ denotes the initial condition.
In Sch\"utz's formula, this would involve a  sum over the positions of the other $N-m$ particles.  The resulting formula is not useful,
and, in particular, not conducive to the $N\to \infty$ limit.

This was overcome by \cite{sasamoto,borFerPrahSasam}, who were able to rewrite the right hand side of \eqref{eqGreen} as a (signed) determinantal point process on a space of Gelfand-Tsetlin patterns.
This allowed them to employ the Eynard-Mehta technology \cite{eynardMehta} to conclude that \eqref{eq:extKernelProb} can be written as the Fredholm determinant $F_t(X_0;a,n)=\det\!\left(I-\bar\chi_a K^\TASEP_t\bar\chi_a\right)$ which (in principle) could be obtained from the operators $\bar\chi_a$ and $K^\TASEP_t$ acting on $\ell^2(\{n_1,\dotsc,n_m\}\times\zz)$ given by
\begin{gather}
\bar\chi_af(n_j,x)=f(n_j,x)\uno{x\leq a_j}\label{eq:chi}\\
\shortintertext{and}
K^\TASEP_t(n_i,x_i;n_j,x_j)=-Q^{n_j-n_i}(x_i,x_j)\uno{n_i<n_j}+\sum_{k=1}^{n_j}\Psi^{n_i}_{n_i-k}(x_i)\Phi^{n_j}_{n_j-k}(x_j),\label{eq:Kt}
\end{gather}
where $Q(x,y)=\frac{1}{2^{x-y}}\uno{x>y}$ and, for $k\leq n-1$,
\begin{equation}\label{eq:defPsi}
\Psi^n_k(x)=\frac1{2\pi\I}\oint_{\Gamma_0}\d w\,\frac{(1-w)^k}{2^{x-X_0(n-k)}w^{x+k+1-X_0(n-k)}}e^{t(w-1)},
\end{equation}
where $\Gamma_0$ is any positively oriented simple loop including the pole at $w=0$ but not the one at $w=1$.
The functions $\Phi_k^{n}(x)$, $k=0,\ldots,n-1$, are defined implicitly by: 
I. The biorthogonality relation $\sum_{x\in\zz}\Psi_k^{n}(x)\Phi_\ell^{n}(x)=\uno{k=\ell}$;
II.  $2^{-x}\Phi^n_k(x)$ is a polynomial of degree at most $n-1$ in $x$ for each $k$.

Except for a few very special choices of initial data, the solution for the $\Phi_k^{n}(x)$ was not discovered until \cite{fixedpt}.
Without discussing their exact form, we state the final result after performing some manipulations to get a nice formula.
Let $\bar Q^{(n)}(x,y)=2^{x-y}\frac{1}{(n-1)!}\prod_{j=0}^{n-1}(x-y-j)$ be the real analytic extension of $Q^n$, let $\nabla^-f(x)=f(x)-f(x-1)$ be the backwards discrete difference operator, and observe that $Q$ is invertible (with inverse $Q^{-1}(x,y)=2\cdot\uno{x=y-1}-\uno{x=y}$).
Define
\begin{align}
 \SM_{-t,-n}(z_1,z_2) &=(e^{-\frac{t}{2}\nabla^-}\tsm Q ^{-n})^*(z_1,z_2)=\frac{1}{2\pi\I} \oint_{\Gamma_0}\d w\, \frac{(1-w)^{n}}{2^{z_2-z_1} w^{n +1 + z_2 - z_1}}e^{t(w-1/2)},\label{def:sm}\\
 \SN_{-t,n} (z_1,z_2) &=\bar Q^{(n)}e^{\frac{t}{2}\nabla^-}(z_1,z_2)=\frac{1}{2 \pi \I} \oint_{\Gamma_{0}} \d w\,\frac{(1-w)^{z_2-z_1 + n - 1}}{2^{z_1-z_2} w^{n}} e^{t(w-1/2)},\label{def:sn}
\end{align}
the first one being defined for all $n\in\nn$ and the second one for $n\geq1$.
Here $(e^{-\frac{t}{2}\nabla^-})_{t\geq0}$ is the semigroup of a Poisson process with jumps to the left at rate $\frac{1}{2}$, which we may think of  as an integral operator on $\ell^2(\zz)$ with kernel $e^{-\frac{t}{2}\nabla^-}\!(x,y)=e^{-\frac{t}{2}}\frac{t^{x-y}}{2^{x-y}(x-y)!}\uno{x\geq y}$; this formula is actually valid for all $t\in\rr$ and it defines the whole group of operators $(e^{-\frac{t}{2}\nabla^-})_{t\in\rr}$ (in particular $e^{-\frac{t}{2}\nabla^-}$ is invertible, with inverse $e^{\frac{t}{2}\nabla^-}$).
Define also, for $n\geq0$,
\begin{equation}\label{eq:sepi23}
{\SN}_{-t,n}^{\epi(X_0)}(z_1,z_2) = \ee_{B_0=z_1}\!\left[ \SN_{-t,n - \tau}(B_{\tau}, z_2)\uno{\tau<n}\right],
\end{equation}
where $\tau$ is defined to be the hitting time of the strict epigraph of the `curve' $\big(X_0(k+1)\big)_{k=0,\dotsc,n-1}$ by a discrete time random walk $B_k$ with transition probabilities $Q$.
Then, as proven in \cite{fixedpt}, the kernel from \eqref{eq:Kt} can be expressed as
\begin{equation}\label{eq:Kt-2}
K^\TASEP_t(n_i,\cdot;n_j,\cdot)=-Q^{n_j-n_i}\uno{n_i<n_j}+(\SM_{-t,-{n}_i})^*{\SN}_{-t,n_j}^{\epi(X_0)}.
\end{equation}

Now let $\cL$ denote the generator of $\TASEP$.
It acts on bounded cylinder (i.e. depending on finitely many coordinates) functions $f\!:\mathcal{W}\longrightarrow\rr$, where $\mathcal{W}=\{X\in\zz^\nn\!:X(1)>X(2)>\dotsc\}$, as follows:
\begin{equation}\label{eq:TASEPgen}
\cL f(X)=\sum_{k\geq1}\uno{X({k-1})-X(k)>1}\tsm\big[f(X(1),\dotsc,X(k)+1,\dotsc)-f(X(1),\dotsc,X(k),\dotsc)\big],
\end{equation}
where we take the convention $X(0)=\infty$ in the indicator $\uno{X(0) - X(1) > 1}$.
Since $f$ is a cylinder function, the sum has only finitely many non-zero terms. 
If $X_t$ is a Markov process with generator $\cL$, its transition probabilities $p(s,X;T,A)=\pp(X_T\in A\mid X_s=X)$ satisfy the \emph{Kolmogorov backward equation}, $\partial_s p(s,X;T,A) = -\cL\, p(s,X;T,A)$ with \emph{final condition} $\lim_{s\uparrow T}  p(s,X;T,A)= \uno{X\in A}$.
The generator $\cL$ acts on the $X$ variable.
In our case, we have a time homogeneous process $ p(s,X;T,A)= \tilde{p}(T-s,X;A)$  and a generating family of sets $A= \{ X(n_j)>a_j, j=1,\ldots,m\}$ and $F_t(X;a,n)= \tilde{p}(t,X;A)$, so the backward equation for $F$ reads,
\begin{align}
\tfrac{\d}{\d t}F&=\cL F\label{eq:bckwd1}
\end{align}
with \emph{initial condition}
\begin{align}
F_0(X;a,n)&=\uno{X(n_1)>a_1,\dotsc X(n_m)>a_m}.\label{eq:bckwd2}
\end{align}
The solution of the Kolmogorov backward equation for a continuous time Markov chain on a countable state space $S$ is unique  under the condition that the supremum over all states of the rate of leaving that state is finite (see, for example, Thm. 2.26 and Cor. 2.34 of \cite{ligg3}, where it is written for the point-to-point probabilities $p_t(X,Y) =\pp_X(X(t)=Y)$; $S$ being countable, the probability measure $p_t(X,A) =P_X(X(t)\in A)=\sum_{Y\in A}p_t(X,Y)$ is therefore prescribed, and satisfies uniquely the backward equation in $t,X$, with initial condition $p_0(X,A)=\uno{X\in A}$).
Our state space $\mathcal{W}$ is \emph{not} countable; however, in TASEP (and PushASEP), the evolution of particles $X(1),\ldots,X(n)$ is unaffected by particles $X(m)$, $m>n$.  
Hence the backward equation \eqref{eq:bckwd1}, \eqref{eq:bckwd2} actually takes place on the countable set $\mathcal{W}_{n_m}=\{(X(1),X(2),\dotsc,X(n_m))\in\zz^\nn\!:X(1)>X(2)>\dotsc>X(n_m)\}$.  
Therefore we have

\begin{theorem}\label{thm:KolmogorovBackwards} 
The unique solution of the \emph{Kolmogorov backward equation} \eqref{eq:bckwd1}, \eqref{eq:bckwd2} for TASEP
is given explicitly by the Fredholm determinant of the kernel $K^\uptext{TASEP}_t$ introduced in \eqref{eq:Kt}:
\begin{equation}\label{detform}
F_t(X_0;a,n)=\det(I-\bar\P_a K^\uptext{TASEP}_t\bar\P_a)_{\ell^2(\{n_1,\dotsc,n_m\}\times\zz)}.
\end{equation}
\end{theorem}

However, a detailed proof along the historical lines sketched above would run about 30 pages.
A natural and important question is whether one could just prove directly that the determinant satisfies the Kolmogorov equation.
The purpose of this article is to provide such a proof, for TASEP and its variant PushASEP.
We comment that the direct proof is relatively short, but far from obvious.
It is not known if it can be obtained directly from the biorthogonal representation \eqref{eq:Kt} without knowledge of the special form of the $\Phi_k^{n}(x)$.
It is also worth noting that there are examples (such as discrete time TASEP, see \cite{mqr-variants}) with explicit $\Phi_k^{n}(x)$ for which the present short proof does not work.

\vskip8pt

TASEP has a rich history, some of which can be found in the introduction of \cite{MR1468391}.
Since that paper there has been significant progress in exact solvability, with much interest stemming from TASEP's role as one of the fundamental models in the KPZ universality class.
There are many results; we mention here only a few besides the ones which fit directly into our story above.
The one-point distribution of TASEP with step initial data ($X_0(i)=-i$, $i\geq1$) was first solved in \cite{johanssonShape} by exploiting a Toeplitz structure which is only available for that choice of initial condition; the main goal of that paper was to prove the now-classic fact that the fluctuations of the position of an appropriately chosen particle converge, as time goes to infinity, to the Tracy-Widom GUE distribution \cite{tracyWidom}.
This was later extended in \cite{johansson} to the multi-point distributions of TASEP, which were shown there to converge to the Airy$_2$ process \cite{prahoferSpohn}.
Much work was devoted during the last fifteen years to extending this type of results to two other choices of initial data: periodic ($X_0(i)=-2i$, $i\in\zz$, which was the subject of the papers \cite{sasamoto,borFerPrahSasam} discussed above) and stationary ($X_0$ corresponding to placing particles on $\zz$ according to a product measure) \cite{imamSasam1,bfs,baikFerrariPeche}.  
In the past ten years there has been a huge effort to extend some of these results to the partially asymmetric simple exclusion process (see, for example, \cite{tracyWidomASEP2,imamSasamDuality,bcs}).
These still depend on very particular initial data.
In the case of TASEP on a ring there is also a huge literature with some recent breakthroughs: \cite{baikLiu-ring}, where asymptotics is done
for a type of step initial data, and \cite{prolhac-spectrum2,prolhac-spectrum1}, which actually computes the entire spectrum.

\section{Kolmogorov equation for TASEP}\label{sec:kolm}

In this section we provide a short self-contained proof of Thm. \ref{thm:KolmogorovBackwards}.
The strategy is to compute the two sides of \eqref{eq:bckwd1} with $F_t$ given by \eqref{detform} and $K^\uptext{TASEP}_t$ given by \eqref{eq:Kt-2}, and check that they are equal.

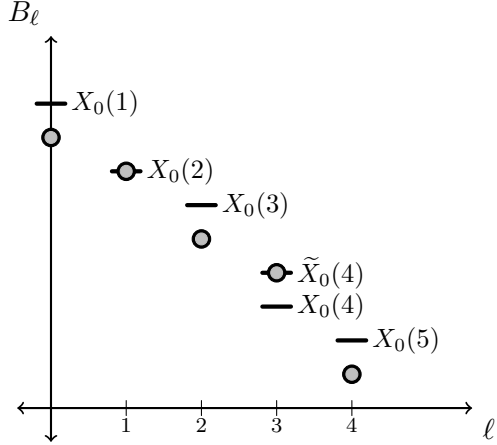
\begin{figure}
\[
\begin{aligned}
\begin{tikzpicture}[scale=1]
\draw [<->,thick] (0,11*0.45) node (yaxis) [above left] {$B_\ell$}
        |- (5.6,0) node (xaxis) [below right] {$\ell$};
\draw [<-,thick] (-0.45,0) -- (0,0);
\draw [->,thick] (0,0) -- (0,-0.45);
\draw[ultra thick,cap=round] (0-0.19,9*0.45)--node (X1) [right] {\small$\,\,X_0(1)$} ++(2*0.19,0);
\draw[ultra thick,cap=round] (1-0.19,7*0.45)--node (X2) [right] {\small$\,\,X_0(2)$} ++(2*0.19,0);
\draw[cap=round] (1,1*0.1) --node (l3) [below] {\scriptsize$1$} (1,-1*0.1);
\draw[ultra thick,cap=round] (2-0.19,6*0.45)--node (X3) [right] {\small$\,\,X_0(3)$} ++(2*0.19,0);
\draw[cap=round] (2,1*0.1) --node (l3) [below] {\scriptsize$2$} (2,-1*0.1);
\draw[ultra thick,cap=round] (3-0.19,3*0.45)--node (X4) [right] {\small$\,\,X_0(4)$} ++(2*0.19,0);
\draw[cap=round] (3,1*0.1) --node (l3) [below] {\scriptsize$3$} (3,-1*0.1);
\draw[ultra thick,cap=round] (3-0.19,4*0.45)--node (tildeX4) [right] {\small$\,\,\wt{X}_0(4)$} ++(2.*0.19,0);
\draw[cap=round] (4,1*0.1) --node (l3) [below] {\scriptsize$4$} (4,-1*0.1);
\draw[ultra thick,cap=round] (4-0.19,2*0.45)--node (X5) [right] {\small$\,\,X_0(5)$} ++(2*0.19,0);

\draw[fill] (0,8*0.45) circle (3.5pt);
\draw[fill,lightgray] (0,8*0.45) circle (2.25pt);
\draw[fill] (1,7*0.45) circle (3.5pt);
\draw[fill,lightgray] (1,7*0.45) circle (2.25pt);
\draw[fill] (2,5*0.45) circle (3.5pt);
\draw[fill,lightgray] (2,5*0.45) circle (2.25pt);
\draw[fill] (3,4*0.45) circle (3.5pt);
\draw[fill,lightgray] (3,4*0.45) circle (2.25pt);
\draw[fill] (4,1*0.45)circle (3.5pt);
\draw[fill,lightgray] (4,1*0.45) circle (2.25pt);
\end{tikzpicture}\\
\end{aligned} 
\]
\caption{A sample path of the random walk $B_\ell$ (circles) and the graph of the initial conditions $\{X_0(i+1)\}_{i=0}^{n-1}$ and $\{\tilde{X}_0(i+1)\}_{i=0}^{n-1}$ which differ only at position $k$ with $\tilde{X}_0(k)=X_0(k)+1$.
The random walk has i.i.d. geometric decrements, and $\tau$ and $\tilde{\tau}$ are the first times the random walk $B_\ell$ hits the strict epigraph of the initial condition curve with $X_0$ and $\tilde{X}_0$ respectively. Note that $\tilde{\tau}=\tau$ unless the hit happens exactly at $\tau=k-1<n$ and $B_{k-1}=\tilde{X}_0(k) = X_0(k)+1$.
In this example $k=4$, $\tau=3$ and $\tilde{\tau}>4$.} \label{fig:example}
\end{figure}

We begin with the right hand side of \eqref{eq:bckwd1}.
We will first consider the effect of moving a single particle, and then we will sum over all particles to obtain the effect of the generator $\cL$. Fix a particle label $k$ and consider the original initial condition $X_0$ and the initial condition $\wt X_0$ where particle $k$ is moved to the right by one, i.e. $\wt X_0(k) \coloneqq X_0(k)+1$ (for now it does not even matter whether or not this particle can actually be moved without violating the strict particle order condition).
We will compare the kernel $K^\uptext{TASEP}_t$ with these two initial conditions; for notational convenience we put tildes on top of all objects when they depend on the modified initial condition (e.g. $\wt K_t^\uptext{TASEP}$ refers to the kernel started from the initial condition $\wt X_0$).
By \eqref{eq:Kt-2} we have
\begin{equation}
\big(\wt K_t^\uptext{TASEP}-K_t^\uptext{TASEP}\big)(n_i,\cdot;n_j,\cdot)=(\SM_{-t,-n_i})^*\big(\SN_{-t,n_j}^{\epi(\wt X_0)}-\SN_{-t,n_j}^{\epi(X_0)}\big),\label{eq:diffKt}
\end{equation}
while by \eqref{eq:sepi23} we have
\begin{equation}
\big(\SN_{-t,n}^{\epi(\wt X_0)}-\SN_{-t,n}^{\epi(X_0)}\big)(z_1,z_2)
=\ee_{B_0=z_1}\!\left[\SN_{-t,n-\tilde\tau}(B_{\tilde\tau},z_2)\uno{\tilde\tau<n}-\SN_{-t,n-\tau}(B_\tau,z_2)\uno{\tau<n}\right].\label{eq:diffSN}
\end{equation}
Recall that $\tau$ means the hitting time of the strict epigraph of $\big(X_0(\ell+1)\big)_{\ell\geq0}$; $\tilde\tau$ is the same but with the modified $\wt X_0$.
It is clear from the definitions that $\tilde\tau=\tau$ unless $\tau=k-1<n$ and $B_{k-1}=X_0(k)+1$ (see Fig. \ref{fig:example} for a visual explanation), and thus the above expectation equals
\begin{equation}
\ee_{B_0=z_1}\!\left[\big(\SN_{-t,n-\tilde\tau}(B_{\tilde\tau},z_2)\uno{\tilde\tau<n}-\SN_{-t,n-\tau}(B_\tau,z_2)\uno{\tau<n}\big)\uno{\tau=k-1,B_{k-1}=X_0(k)+1}\right]\!\uno{k\leq n},
\end{equation}
which by the Markov property we may write as $f_k(z_1)g^{(n)}_{k}(z_2)$ with
\begin{gather}
f_k(z)=\pp_{B_0=z}\!\left(\tau=k-1,\,B_{k-1}=X_0(k)+1\right)\label{eq:fk}\\
\shortintertext{and}
g^{(n)}_{k}(z)=\Big(\ee_{B_{k-1}=X_0(k)+1}\!\left[\SN_{-t,n-\tilde\tau^{(k)}}(B_{\tilde\tau^{(k)}},z)\uno{\tilde\tau^{(k)}<n}\right]-\SN_{-t,n-k+1}(X_0(k)+1,z)\Big)\uno{k\leq n},~\quad\label{eq:gk}
\end{gather}
where the subscript in expectation on the right hand side indicates that the random walk is started now at $X_0(k)+1$ at time $k-1$ and the superscript in $\tilde\tau^{(k)}$ indicates that the hitting time is now restricted to times in $\{k,k+1,\dotsc\}$ (we will use the same notation for other hitting times below).
This means that 
\begin{equation}
\Delta^{\!(k)}_{\tts i,j}\coloneqq\big(\wt K_t^\uptext{TASEP}-K_t^\uptext{TASEP}\big)(n_i,\cdot;n_j,\cdot)=(\SM_{-t,-n_i})^*f_k\otimes g^{(n_j)}_{t,k}.\label{eq:KcompTASEP}
\end{equation}
Crucially, this is a rank one kernel (acting on $\ell^2(\zz)$), and then the same is true of the extended kernel $\wt K^\uptext{TASEP}_t-K^\uptext{TASEP}_t$ when thought of as acting on the extended space $\ell^2(\{n_1,\dotsc,n_m\}\times\zz)$, so using the fact that $\det(I+A+B)-\det(I+A)=\det(I+A)\tr[(I+A)^{-1}B]$ for $B$ any rank one operator, we deduce that
\begin{multline}
\det(I-\bar\P_a\wt K^\uptext{TASEP}_t\bar\P_a)-\det(I-\bar\P_aK^\uptext{TASEP}_t\bar\P_a)\\
=-\det(I-\bar\P_aK^\uptext{TASEP}_t\bar\P_a)\tr[(I-\bar\P_aK^\uptext{TASEP}_t\bar\P_a)^{-1}\bar\P_a\Delta^{\!(k)}\bar\P_a]
\end{multline}
(here $\Delta^{\!(k)}=\wt K_t^\uptext{TASEP}-K_t^\uptext{TASEP}$ as in \eqref{eq:KcompTASEP}).
This difference corresponds to the bracket in the $k$-th term of the sum \eqref{eq:TASEPgen} defining the action of $\cL$ on our function $F_t(X_0;a,n)$.
The conclusion then is that if $J\subseteq\nn$ is the set of labels of particles which can be moved (i.e. those so that $X_0(k-1)-X_0(k)>1$), then $\cL F_t(X_0;a,n)$ equals the sum over $k\in J$ of the above difference, i.e.
\begin{equation}
\cL F_t(X_0;a,n)
=-\det(I-\bar\P_aK^\uptext{TASEP}_t\bar\P_a)\tr\!\Big[(I-\bar\P_aK^\uptext{TASEP}_t\bar\P_a)^{-1}\bar\P_a\big({\textstyle\sum_{k\in J}}\Delta^{\!(k)}\big)\bar\P_a\Big].
\end{equation}
Note that, as above, the sum only has finitely many non-zero terms, because $\Delta^{\!(k)}=0$ for all $k>n_m$ (as is clear, for instance, from \eqref{eq:gk} and \eqref{eq:KcompTASEP}).
Finally we note if $k\notin J$ (i.e. if $X_0(k)+1=X_0(k-1)$), then since the random walk $B_\cdot$ has increments $B_{t-1} - B_{t} \geq 1$, we have $\{B_{k-1}=X_0(k)+1\}\subseteq\{B_{k-2}\geq X_0(k-1)+1\}\subseteq\{\tau\leq k-2\}$ and thus $f_k(z)=0$. Thus summing over $j \notin J$ has no effect, and we have then
\begin{equation}\label{eq:RHSfinal}
\cL F_t(X_0;a,n)
=-\det(I-\bar\P_aK^\uptext{TASEP}_t\bar\P_a)\tr\!\Big[(I-\bar\P_aK^\uptext{TASEP}_t\bar\P_a)^{-1}\bar\P_a\big({\textstyle\sum_{k\geq1}}\Delta^{\!(k)}\big)\bar\P_a\Big].
\end{equation}

Now we turn to the left hand side of \eqref{eq:bckwd1}.
In general, if a kernel $K_h$ depends smoothly on a parameter $h$ one has $\frac{\d}{\d h}\det(I-K_h)=-\det(I-K_h)\tr[(I-K_h)^{-1}\frac{\d}{\d h}K_h]$.
Thus
\[\tfrac{\d}{\d t}F_t(X_0;a,n)=-\det(I-\bar\P_aK^\uptext{TASEP}_t\bar\P_a)\tr\!\Big[(I-\bar\P_aK^\uptext{TASEP}_t\bar\P_a)^{-1}\bar\P_a\big(\tfrac{\d}{\d t}K^\uptext{TASEP}_t\big)\bar\P_a\Big]\]
and then in view of \eqref{eq:RHSfinal} all we need in order to deduce \eqref{eq:bckwd1} is to check that
\begin{equation}\label{eq:toprove}
\tfrac{\d}{\d t}K^\uptext{TASEP}_t={\textstyle\sum_{k\geq1}}\Delta^{\!(k)}.
\end{equation}
We have from \eqref{eq:Kt-2} that
\[\Lambda_{i,j}\coloneqq\tfrac{\d}{\d t}K^\uptext{TASEP}_t(n_i,\cdot;n_j,\cdot)=\big(\tfrac{\d}{\d t}(\SM_{-t,-n_i})^*\big)\SN_{-t,n_j}^{\epi(X_0)}+(\SM_{-t,-n_i})^*(\tfrac{\d}{\d t}\SN_{-t,n}^{\epi(X_0)}).\]
From \eqref{def:sm} we get
\begin{align}
\tfrac{\d}{\d t}(\SM_{-t,-n})^*&=-\tfrac12\nabla^-(\SM_{-t,-n}\big)^*=-\tfrac12(\SM_{-t,-n}\big)^*\nabla^-\\
\shortintertext{and from \eqref{def:sn} we get $\tfrac{\d}{\d t}\SN_{-t,n}=\tfrac12\SN_{-t,n}\nabla^-=\tfrac12\nabla^-\SN_{-t,n}$, which means}
\tfrac{\d}{\d t}\SN_{-t,n}^{\epi(X_0)}(z_1,z_2)&=\ee_{B_0=z_1}\tsm\big[\tfrac12\nabla^-\SN_{-t,n-\tau}(B_\tau,z_2)\uno{\tau<n}\big]
\end{align}
(we have used here that $\nabla^-$ commutes with $(\SM_{-t,-n})^*$ and $\SN_{-t,n}$ since they all have Toeplitz kernels).
This gives
\[\Lambda_{i,j}=\tfrac12(\SM_{-t,-n_i})^*H_{n_j}\]
with
\begin{equation}\label{eq:DeltaTASEP}
\begin{aligned}
H_{n}(z_1,z_2)&=-\nabla^-\SN_{-t,n}^{\epi(X_0)}(z_1,z_2)+\ee_{B_0=z_1}\tsm\big[\nabla^-\SN_{-t,n-\tau}(B_\tau,z_2)\uno{\tau<n}\big]\\
&=\SN_{-t,n}^{\epi(X_0)}(z_1-1,z_2)-\ee_{B_0=z_1}\tsm\big[\SN_{-t,n-\tau}(B_\tau-1,z_2)\uno{\tau<n}\big].
\end{aligned}
\end{equation}
Shifting the curve defined by $X_0$ up by one, the first term on the second line can be expressed as $\ee_{B_0=z_1}\tsm\big[\SN_{-t,n-\hat\tau}(B_{\hat\tau}-1,z_2)\uno{\hat\tau<n}\big]$ with $\hat\tau$ the hitting time of the strict epigraph of $\big(X_0(m+1)+1\big)_{m\geq0}$, and thus
\[H_n(z_1,z_2)=\ee_{B_0=z_1}\tsm\big[\SN_{-t,n-\hat\tau}(B_{\hat\tau}-1,z_2)\uno{\hat\tau<n}-\SN_{-t,n-\tau}(B_{\tau}-1,z_2)\uno{\tau<n}\big].\]
Note that $\hat\tau\geq\tau$, so the difference inside the brackets vanishes for $\tau\geq n$, while if $\tau<n$ and $B_\tau\geq X_0(\tau+1)+2$ we have $\hat\tau=\tau$, so the difference vanishes again.
We deduce that
\begin{equation}
\begin{split}
 H_n(z_1,z_2)&=\sum_{k=0}^{n-1}\pp_{B_0=z_1}\!\left(\tau=k,\,B_{k}=X_0(k+1)+1\right)\\
&\qquad\times\left[\ee_{B_{k}=X_0(k+1)+1}\!\left[\SN_{-t,n-\hat\tau^{(k)}}(B_{\hat\tau^{(k)}}-1,z_2)\uno{\hat\tau^{(k)}<n}\right]-\SN_{-t,n-k}(X_0(k+1),z_2)\right]\\
&=\sum_{k=0}^{n-1}f_{k+1}(z_1)\ts\widehat g^{(n)}_{k+1}\tsm(z_2)
\end{split}
\end{equation}
with $f_k$ as in \eqref{eq:fk} and $\widehat g^{(n)}_{k+1}\tsm(z)$ given by the bracket in the middle line (see also the comment after \eqref{eq:gk}).

Shifting now the curve defining $\hat\tau$ back to $X_0$ in the last expectation we get
\begin{equation}\label{eq:wgk}
\widehat g^{(n)}_{k}\tsm(z)=\ee_{B_{k-1}=X_0(k)}\!\left[\SN_{-t,n-\tau}(B_{\tau},z)\uno{\tau<n}\right]-\SN_{-t,n-k+1}(X_0(k),z).
\end{equation}
We have then
\begin{equation}
\tfrac{\d}{\d t}K^\uptext{TASEP}_t(n_i,\cdot;n_j,\cdot)=\tfrac12{\textstyle\sum_{k=1}^{n_j}}(\SM_{-t,-n_i})^*f_k(z_1)\otimes\widehat g^{(n_j)}_{k}\tsm(z_2).\label{eq:Kderiv}
\end{equation}
Comparing this formula and \eqref{eq:KcompTASEP}, and since $g^{(n)}_k\equiv0$ for $k\geq n$ by its definition \eqref{eq:gk}, all that remains in order to prove \eqref{eq:toprove} is that
\begin{equation}
\widehat g^{(n)}_k=2g^{(n)}_k\quad\uptext{ for all $k\leq n$}.\label{eq:finaltoprove}
\end{equation}
But from \eqref{def:sn} we have $\SN_{-t,n}e^{-\frac{t}{2}\nabla^-}=\SN_{0,n}=\bar Q^{(n)}$, so applying $e^{-\frac{t}{2}\nabla^-}$ on the right in \eqref{eq:finaltoprove} and using the definitions \eqref{eq:gk} and \eqref{eq:wgk} of $g^{(n)}_k$ and $\widehat g^{(n)}_k$ we see that \eqref{eq:finaltoprove} is equivalent to
\begin{multline}\label{eq:noRt}
\ee_{B_{k-1}=X_0(k)}\!\left[\bar Q^{(n-\tau^{(k)})}(B_{\tau^{(k)}},z)\uno{\tau^{(k)}<n}\right]-\bar Q^{(n-k+1)}(X_0(k),z)\\
=2\tts\ee_{B_{k-1}=X_0(k)+1}\!\left[\bar Q^{(n-\tilde\tau^{(k)})}(B_{\tilde\tau^{(k)}},z)\uno{\tilde\tau^{(k)}<n}\right]-2\tts\bar Q^{(n-k+1)}(X_0(k)+1,z)
\end{multline}
for $k\leq n$.
By definition of $\bar Q^{(n)}$, if we multiply the equation by $2^{-z}$ then both sides are polynomials in $z$, so it is enough to prove the equality for $z<X_0(n)$.
But it is easy to check that 
\begin{equation}\label{eq:QQbar}
\bar Q^{(n)}(x,y)=Q^n(x,y)\qquad\uptext{for all }x-y\geq1.
\end{equation}
(see \fpeqref{eq:QQext}), so for such $z$, since $X_0(i)$ is decreasing in $i$ then $z$ is also smaller than $B_\tau$, $B_{\tilde\tau}$ and $X_0(k)$ in \eqref{eq:noRt}, and we deduce that all the $\bar Q^{(\ell)}$'s in the identity can be replaced by $Q^{\ell}$.
As a consequence, what we need to prove is that
\begin{multline}
\ee_{B_{k-1}=X_0(k)}\!\left[Q^{n-\tau^{(k)}}(B_{\tau^{(k)}},z)\uno{\tau^{(k)}<n}\right]-Q^{n-k+1}(X_0(k),z)\\
=2\tts\ee_{B_{k-1}=X_0(k)+1}\!\left[Q^{n-\tilde\tau^{(k)}}(B_{\tilde\tau^{(k)}},z)\uno{\tilde\tau^{(k)}<n}\right]-2\tts Q^{n-k+1}(X_0(k)+1,z),
\end{multline}
or simply
\[\pp_{B_{k-1}=X_0(k)}\!\left(B_n=z,\,\tau^{(k)}\geq n\right)=2\pp_{B_{k-1}=X_0(k)+1}\!\left(B_n=z,\,\tilde\tau^{(k)}\geq n\right),\]
which is easy to see since the walk takes Geom$[\frac12]$ steps to the left: in fact, since $X_0(k+1)<X_0(k)$ and $\tau^{(k)}$ and $\tilde\tau^{(k)}$ may only differ on the event that the walk hits at time $k-1$, then for $k<n$, decomposing according to the first step, 
\begin{align}
&\textstyle\pp_{B_{k-1}=X_0(k)+1}\!\left(B_n=z,\,\tilde\tau^{(k)}\geq n\right)=\sum_{y\leq X_0(k+1)}\tsm Q(X_0(k)+1,y)\tts\pp_{B_{k}=y}\!\left(B_n=z,\,\tilde\tau^{(k)}\geq n\right)\\
&\hspace{0.05in}=\textstyle\sum_{y\leq X_0(k+1)}\tsm\tfrac12Q(X_0(k),y)\tts\pp_{B_{k}=y}\!\left(B_n=z,\,\tau^{(k)}\geq n\right)
=\tfrac12\pp_{B_{k-1}=X_0(k)}\!\left(B_n=z,\,\tau^{(k)}\geq n\right);
\end{align}
the case $k=n$ is even simpler.
This yields the desired equality \eqref{eq:finaltoprove}.

\section{Initial condition}

Now we check that the initial condition \eqref{eq:bckwd2} is satisfied.
It is simpler to deal first with the one-point case ($m=1$), so we start there.

\subsection{One-point case}

Let
\[K_t^{(n)}=K^\uptext{TASEP}_t(n,\cdot;n,\cdot).\]
We need to prove that
\begin{equation}\label{eq:init-toshow}
  \det\!\left(I-\bar\chi_aK_0^{(n)}\bar\chi_a\right)_{\ell^2(\zz)}=\uno{X_0(n)>a}.
\end{equation}

\begin{rem}
In the one point case, the initial condition \eqref{eq:init-toshow} can be checked rather easily from the biorthogonal ensemble representation \eqref{eq:Kt}, as mentioned by one of the referees; the proof using this biorthogonalization is given below.
However, since our goal in this paper is to prove \eqref{detform} directly from the representation \eqref{eq:Kt-2} (and thus avoid having to use the result of \cite{fixedpt} to connect this formula back to \eqref{eq:Kt}), we provide the direct (but relatively more complicated) proof from \eqref{eq:Kt-2} in this subsection.

\noindent From \eqref{eq:KcompTASEP} we have $K^{(n)}_0(z_1,z_2)=\sum_{k=0}^{n-1}\Psi^n_k(z_1)\Phi^n_k(z_2)$ with $t=0$ in the definition \eqref{eq:defPsi}.
Let $\Psi,\Phi\!:\zz\times\{0,\dots,n-1\}\longrightarrow\rr$ be given by $\Psi(z,k)=\Psi^n_k(z)$ and $\Phi(z,k)=\Phi^n_k$, so that the left hand side of \eqref{eq:init-toshow} equals 
\[\det\!\left(I-\bar\chi_a\Psi\Phi^*\bar\chi_a\right)_{\ell^2(\zz)}=\det(I-\Phi^*\bar\P_a\Psi),\]
where we have used the cyclic property of the Fredholm determinant.
Note that on the determinant on the right hand side is finite dimensional, with the $n\times n$ matrix inside being given by 
\[\Phi^*\bar\P_a\Psi(\ell,k)=\sum_{z\leq a}\Psi^n_k(z)\Phi^n_\ell(z),\qquad \ell,k=0,\dotsc,n-1.\]
Assume first that $a<X_0(n)$.
Then $a<X_0(n-k)-k$ for all $k=0,\dotsc,n$, which implies from \eqref{eq:defPsi} that $\Psi^n_k(z)=0$ for all $z\leq a$ and all such $k$, so $\det(I-\Phi^*\bar\P_a\Psi)=\det(I)=1$ as desired.
Assume otherwise that $a\geq X_0(n)$.
From \eqref{eq:defPsi} it is straightforward to check that $\Psi^n_0(z)=\uno{z=X_0(n)}$ and thus the first column of $\Phi^*\bar\P_a\Psi$ is given by $\Phi^*\bar\P_a\Psi(\ell,0)=\sum_{z\leq a}\Psi^n_k(z)\Phi^n_\ell(z)=\Phi^n_\ell(X_0(n))$.
But the biorthogonalization condition $\sum_{z\in\zz}\Psi^n_0(z)\Phi^n_\ell(z)=\uno{\ell=0}$ implies that $\Phi^n_\ell(X_0(n))=\uno{\ell=0}$.
We deduce that the first column of $I-\Phi^*\bar\P_a\Psi$ equals zero, so $\det(I-\Phi^*\bar\P_a\Psi)$ vanishes as desired.
\end{rem}

We turn now to the proof of \eqref{eq:init-toshow} directly from \eqref{eq:Kt-2}.
Consider first the case $a\leq X_0(n)-1$, which implies that in \eqref{eq:init-toshow} we need to evaluate $K^{(n)}_0(z_1,z_2)$ only for $z_1\leq
 X_0(n)-1$.
By definition of the hitting time $\tau$ in \eqref{eq:sepi23} we have $\SN^{\epi(X_0)}_{0,n}(z_1,z_2)=\bar Q^{(n)}(z_1,z_2)$ for $z_1>X_0(1)$, so we may write
\begin{equation}
K_0^{(n)}=Q^{-n}\chi_{X_0(1)}\bar Q^{(n)}+Q^{-n}\bar\chi_{X_0(1)}\SN_{0,n}^{\epi(X_0)}.\label{eq:K0}
\end{equation}
Now
\begin{equation}\label{eq:Qinv}
Q^{-n}(x,y)=(-1)^{y-x+n}\tts2^{y-x}\ttsm{\textstyle\binom{n}{y-x}}\uno{0\leq y-x\leq n}
\end{equation} 
so, in particular, $Q^{-n}(x,y)=0$ for $y-x>n$.
If $y>X_0(1)$ and $z_1\leq X_0(n)-1$ we have $y-z_1>X_0(1)-X_0(n)+1\geq n$, so $Q^{-n}(z_1,y)=0$ and the first term on the right hand side of \eqref{eq:K0} vanishes.
For the second term we consider $\SN^{\epi(X_0)}_{0,n}(y,z_2)=\sum_{k=0}^{n-1}\ee_{B_0=y}\big[\SN_{0,n-k}(B_k,z_2)\uno{\tau=k}\big]$ with $y\leq X_0(1)$ and note that $\tau=k$ implies $B_k>X_0(k+1)$, which cannot happen unless the starting location satisfies $y>X_0(k+1)+k$ (because the walk takes downward steps of size at least $1$).
But for such $y$ (and $z_1\leq X_0(n)-1$) we have $y-z_1>X_0(k+1)+k-X_0(n)+1\geq n$, and thus again we get $Q^{-n}(z_1,y)=0$.
We have proved that
\begin{equation}
\bar\chi_a K^{(n)}_0=0\qquad\uptext{when $a<X_0(n)$},\label{eq:easyini}
\end{equation}
and then the left hand side of \eqref{eq:init-toshow} equals $\det(I)=1$ for $a<X_0(n)$, as desired.

Next we turn to the case $a\geq X_0(n)$, and focus on the column corresponding to the index $X_0(n)$ in $\bar\P_aK^{(n)}_0\bar\P_a$, that is $\big(K^{(n)}_0(z,X_0(n))\big)_{z\leq a}$.
We claim that, in this case,
\begin{equation}
K^{(n)}_0(z,X_0(n))=\uno{z=X_0(n)}\qquad\text{for all $z\in\zz$}.\label{eq:zleqa}
\end{equation}
If that is the case then the matrix $I-\bar\P_aK^{(n)}_0\bar\P_a$ has a column which is identically 0, and thus its determinant vanishes as required.

In order to prove \eqref{eq:zleqa}, observe first that the argument we used to prove \eqref{eq:easyini} actually shows that $K^{(n)}_0(z,X_0(n))=0$ for all $z<X_0(n)$.
Consider next the case $z=X_0(n)$.
For the first term in \eqref{eq:K0}, a simple residue computation using the contour integral formulas for $Q^{-n}$ and $\bar Q^{(n)}$ which follow from \eqref{def:sm} and \eqref{def:sn} with $t=0$ gives (remembering $X_0(n)-X_0(1)\leq 1-n$)
\begin{equation}\label{eq:init1}
\begin{split}
Q^{-n}\P_{X_0(1)}\bar Q^{(n)}(X_0(n),X_0(n))&\textstyle=\itwopii{2}\oint \d w\oint \d v\,\frac{(1-w)^n(1-v)^{X_0(n)-X_0(1)+n-1}}{w^{X_0(n)-X_0(1)+n}v^n}\frac{1}{1-v-w}\\
&=\uno{X_0(n)-X_0(1)=1-n}.
\end{split}
\end{equation}
Since we want to show that $K^{(n)}_0(X_0(n),X_0(n))=1$, then from \eqref{eq:init1} together with \eqref{eq:K0} we see that we need to show that
\begin{equation}
Q^{-n}\bar\chi_{X_0(1)}\SN_{0,n}^{\epi(X_0)}(X_0(n),X_0(n))=\uno{X_0(n)-X_0(1)<1-n}.\label{eq:QnbPSN}
\end{equation}
We consider the two possibilities for $X_0(n)-X_0(1)$ separately.
Assume first that $X_0(1)-X_0(n)=n-1$.
This means that the first $n$ particles are packed one next to the other, i.e. $X_0(k)=X_0(1)-k+1$, and thus, since the walk starts below $X_0(1)$ and takes downward steps of size at least $1$, it cannot hit the strict epigraph of the curve by time $n$, so $\bar\P_{X_0(1)}\SN_{0,n}^{\epi(X_0)}=0$ as desired.
The other possibility is that $X_0(1)-X_0(n)>n-1$.
Using $\SN_{0,n}=\bar Q^{(n)}$ we may express the kernel $\SN^{\epi(X_0)}_{0,n}(y,X_0(n))$ as $\sum_{k=0}^{n-1}\ee_{B_0=y}\tsm\big[\bar{Q}^{(n-k)}(B_k,X_0(n))\uno{\tau=k}\big]$ and then note that, since inside the expectation $B_k$ is the location where the strict epigraph is hit, $B_k-X_0(n)\geq X_0(k+1)+1-X_0(n)\geq n-k\geq1$, so by \eqref{eq:QQbar} we may replace $\bar{Q}^{(n-k)}$ by $Q^{n-k}$ there to get $\sum_{k=0}^{n-1}\pp_{B_0=y}(B_n=X_0(n),\,\tau=k)$.
But $B_n=X_0(n)$ implies $B_{n-1}\geq X_0(n)+1$, which in turn implies $\tau\leq n-1$, so the last sum equals simply $Q^n(y,X_0(n))$, and then 
\[\textstyle Q^{-n}\bar\P_{X_0(1)}\SN^{\epi(X_0)}_{0,n}(X_0(n),X_0(n))=\sum_{y\leq X_0(1)}Q^{-n}(X_0(n),y)Q^n(y,X_0(n)).\]
And the last sum can be extended to all $y$'s because if $y>X_0(1)$ then, under our assumption $X_0(n)-X_0(1)<1-n$ and using \eqref{eq:Qinv}, the first factor in the summand vanishes.
This shows that, in this case, $Q^{-n}\bar\P_{X_0(1)}\SN^{\epi(X_0)}_{0,n}(X_0(n),X_0(n))=1$ as desired.

Up to here we have proved \eqref{eq:zleqa} for $z\leq X_0(n)$.
Consider finally the case $z>X_0(n)$.
The exact same argument used in the last paragraph to prove \eqref{eq:QnbPSN} also shows 
\[Q^{-n}\bar\P_{X_0(1)}\SN_{0,n}^{\epi(X_0)}(z,X_0(n))=Q^{-n}\chi_{X_0(1)}\bar Q^{(n)}(z,X_0(n))\uno{X_0(1)-X_0(n)>n-1}.\]
In particular, in the case $X_0(1)-X_0(n)>n-1$ we get using \eqref{eq:K0} that $K^{(n)}_0(z,X_0(n))=Q^{-n}\bar Q^{(n)}(z,X_0(n))$ which equals $0$ because in fact $Q^{-n}\bar Q^{(n)}\equiv0$ (see \fpeqref{eq:QinvQext}).
It remains to prove that the first term in \eqref{eq:K0} vanishes when $X_0(1)-X_0(n)=n-1$, but this is straightforward, because in this case we have (using the contour integral formulas \eqref{def:sm} and \eqref{def:sn} as before)
\[\textstyle Q^{-n}\P_{X_0(1)}\SN_{0,n}(z,X_0(n))=\itwopii{2}\oint \d w\oint \d v\,\frac{(1-w)^n}{w^{z-X_0(1)+n}v^n}\frac{1}{1-v-w}=\inv{2\pi\I}\oint \d w\,\frac{1}{w^{z-X_0(1)+n}}=0\]
since $z-X_0(1)+n>X_0(n)-X_0(1)+n=1$.

\subsection{Multi-point case}\label{sec:multpt}

We let $L_{i,j}=K^\uptext{TASEP}_0(n_i,\cdot;n_j,\cdot)$ and think of $K^\uptext{TASEP}_t$ as an operator-valued matrix $L$ with entries $L_{i,j}$.
The initial condition \eqref{eq:bckwd2} will follow from the arguments in the one-point case and the identity
\begin{equation}\label{eq:extending}
L_{i,j}=-Q^{n_j-n_i}\uno{n_i<n_j}+Q^{n_j-n_i}K^{(n_i)}_t,
\end{equation}
which is (\fprefst{eq:checkQKK}) in \cite{fixedpt}.
Consider first the case $a_i\leq X_0(n_i)-1$ for each $i$.
Here \eqref{eq:easyini} gives $\bar\P_{a_i}K^{(n_i)}=0$, and then in the sum $\bar\P_{a_i}Q^{n_j-n_i}K^{(n_i)}_t(z_1,z_2)=\uno{z_1\leq a_i}\sum_{y\in\zz}Q^{n_j-n_i}(z_1,y)K^{(n_i)}_0(y,z_2)$ only terms with $y>a_i$ survive.
On the other hand, we have $Q^{m}(z,y)=0$ for $z<y$ (any $m\in\zz$) so $Q^{n_j-n_i}(z_1,y)=0$ for $z_1\leq a_i<y$.
This shows that $\bar\P_{a_i}L_{i,j}\bar\P_{a_j}=-Q^{n_j-n_i}\uno{n_i<n_j}$, and since this holds for all $i$, we get $F_0(X_0;a,n)=\det(I-\bar\P_aL\bar\P_a)$ with $L$ a strictly upper triangular operator-valued matrix, which yields $F_0(X_0;a,n)=1$ as desired.

Suppose next that $a_j\geq X_0(n_j)$ for some $j$, so we need to prove $\det(I-\bar\P_aL\bar\P_a)=0$.
We may assume without loss of generality that $n_j$ is the largest particle label among $\{n_1,\dotsc,n_m\}$ for which this inequality holds.
We will focus on the $j$-th column of $\bar\P_aL\bar\P_a$, and more precisely on the subcolumn of this column of kernels corresponding to the index $X_0(n_j)$, that is $v^{(j)}\coloneqq\big(v^{(j)}_\ell(z)\big)_{z\in\zz,\,i=1,\dotsc,m}$ with $v^{(j)}_\ell(z)=\bar\P_{a_i}L_{i,j}(z,X_0(n_j))$.
In this case \eqref{eq:zleqa} gives us $K_0^{(n_j)}(z,X_0(n_j))=\uno{z=X_0(n_j)}$ for all $z$, and then $Q^{n_j-n_i}K_0^{(n_j)}(z,X_0(n_j))=Q^{n_j-n_i}(z,X_0(n_j))$, so that using \eqref{eq:extending} we get 
\[v^{(j)}_i(z)=Q^{n_j-n_i}(z,X_0(n_j))\uno{n_i\geq n_j}\uno{z\leq a_i}.\]
The case $i=j$ is direct and yields $v^{(j)}_j(z)=\uno{z=X_0(n_j)}$.
Otherwise, if $n_i>n_j$, then necessarily $a_i<X_0(n_i)$ by our choice of $j$, and hence for $z\leq a_i$ we have $X_0(n_j)-z>X_0(n_j)-X_0(n_i)\geq n_i-n_j-1$ so as before $v^{(j)}_i(z)=0$.
The conclusion is that $v^{(j)}_i(z)=\uno{i=j}\uno{z=X_0(n_j)}$.
But then $I-\bar\P_aL\bar\P_a$ has a column (namely $v^{(j)}$) which is identically $0$, and thus $F_0(X_0;a,n)=\det(I-\bar\P_aL\bar\P_a)=0$ as required.

\section{PushASEP}

We consider now PushASEP, a generalization of TASEP and the Toom model \cite{DLSS91} introduced in \cite{bp-push}.
Again we have particles on the lattice $\zz$ with positions $X_t(1)>X_t(2)>\dotsm$.
Each particle attempts to jump one step to the right at rate $r$, with jumps being permitted only if the neighboring site is empty; this is the same TASEP dynamics considered above (except run at rate $r$).
On the other hand, each particle has another (independent) exponential clock running at rate $\ell$, and when it rings the particle jumps to the nearest vacant site on its left.
Relabeling the particles in order to keep the ordering after such a jump, we may think of the effect of a particle jumping left as pushing all its left neighbors one step to the left.
TASEP and the Toom model (or PushTASEP) are recovered from PushASEP by setting $\ell=0$ and $r=0$, respectively.

The generator $\cL_{r,\ell}$ of this process can be written similary to the one for TASEP: letting $b(k)$ denote the length of the block of nearest neighbor occupied sites lying to the left of (and including) particle $k$, meaning that $X_0(k+j)=X_0(k)-j$ for $j=0,\dotsc,b(k)-1$ and $X_0(m+b(k))<X_0(k)-b(k)$,
\begin{multline}\label{eq:pushASEPgen}
\cL_{r,\ell} f(X)=r\sum_{k\geq1}\uno{X({k-1})-X(k)>1}\tsm\big[f(\dotsc,X(k)+1,\dotsc)-f(\dotsc,X(k),\dotsc)\big]\\
+\ell\sum_{k\geq1}\big[f(\dotsc,X(k)-1,\dotsc,X({k+b(k)})-1,\dotsc)-f(\dotsc,X(k),\dotsc,X({k+b(k)}),\dotsc)\big],
\end{multline}
where again we take $X(0)=\infty$.
As for the TASEP case we define
\begin{equation}
F_t(X;a,n)=\pp_{X}\big(X_t(n_1)>a_1,\dotsc,X_t(n_m)>a_m\big);
\end{equation}
again we have that $F_t(\cdot;a)$ satisfies the \emph{Kolmogorov backward equation}
\noeqref{eq:bckwdpush1}\noeqref{eq:bckwdpush2}
\begin{align}
\tfrac{\d}{\d t}F&=\cL_{r,\ell} F,\label{eq:bckwdpush1}\\
\shortintertext{with initial condition}
F_0(X;a,n)&=\uno{X(1)>a_1,\dotsc X(n)>a_n}.\label{eq:bckwdpush2}
\end{align}

\begin{thm}\label{thm:PushASEPbckwdsoln}
The unique solution of \eqref{eq:bckwdpush1}, \eqref{eq:bckwdpush2} is given by
\[F_t(X_0;a,n)=\det(I-\bar\P_a K^\uptext{PushASEP}_t\bar\P_a)_{\ell^2(\{n_1,\dotsc,n_m\}\times\zz)}\]
with 
\begin{equation}\label{eq:Kt-2-alt}
K^\uptext{PushASEP}_t(n_i,\cdot;n_j,\cdot)=-Q^{n_j-n_i}\uno{n_i<n_j}+(\SM_{-t,-{n}_i})^*{\SN}_{-t,n_j}^{\epi(X_0)},
\end{equation}
where the operators $\SM_{-t,-n}$ and $\SN_{-t,n}$ are now given by
\begin{align}
 \SM_{-t,-n}(z_1,z_2) \label{def:sm-push}
 &= (e^{-\frac12rt\nabla^-+2\ell t\nabla^+}Q^{-n})^*(z_1,z_2)\\
 &=\frac{1}{2\pi\I} \oint_{\Gamma_0}\d w\,\frac{(1-w)^{n}}{2^{z_2-z_1} w^{n +1 + z_2 - z_1}}e^{t[r(w-1/2)+\ell(1/w-2)]},\\
 \SN_{-t,n} (z_1,z_2)
 &=\bar Q^{(n)}e^{\frac12rt\nabla^--2\ell t\nabla^+} (z_1,z_2) \label{def:sn-push}\\
 &= \frac{1}{2 \pi \I} \oint_{\Gamma_{0}} \d w\,\frac{(1-w)^{z_2-z_1 + n - 1}}{2^{z_1-z_2} w^{n}} e^{t[r(w-1/2)+\ell(2-1/(1-w))]};
\end{align}
here $\nabla^+f(x)=f(x+1)-f(x)$ and $e^{-\frac12rt\nabla^-+2\ell t\nabla^+}$ is defined similarly to $e^{-\frac{t}{2}\nabla^-}$ after \eqref{def:sn}.
\end{thm}

This result, which is new, can be proved along the same lines as the solution of TASEP given in \cite{fixedpt}; such a proof will appear in \cite{mqr-variants} (where convergence to the KPZ fixed point will also be proved).
Our goal here, as for the Kolmorogov equation for TASEP in Thm. \ref{thm:KolmogorovBackwards}, is to prove directly that the Fredholm determinant provides a solution to the Kolmogorov backward equation \eqref{eq:bckwdpush1}, \eqref{eq:bckwdpush2}.

Note first that for $t=0$ the formulas for $\SM_{-t,-n}$ and $\SN_{-t,n}$ given in Thm. \ref{thm:PushASEPbckwdsoln} are the same as those for TASEP, so the initial condition \eqref{eq:bckwdpush2} follows from the TASEP proof.
Next we claim that in order to prove that the Fredholm determinant solves \eqref{eq:bckwdpush1} it is enough to consider the case of PushTASEP, i.e. $r=0$, since the TASEP part of the dynamics was already handled in Sec. \ref{sec:kolm}.
To see this it is convenient to define $F_{t_1,t_2}(X_0;a)=\det(I-\bar\P_a K^\uptext{PushASEP}_{t_1,t_2}\bar\P_a)_{\ell^2(\{n_1,\dotsc,n_m\}\times\zz)}$ with
\[K^\uptext{PushASEP}_{t_1,t_2}(n_i,\cdot;n_j,\cdot)=-Q^{n_j-n_i}\uno{n_i<n_j}+(\SM_{-t_1,-t_2,-{n}_i})^*{\SN}_{-t_1,-t_2,n_j}^{\epi(X_0)},\]
$\SM_{-t_1,-t_2,-n}=(e^{-\frac12rt_1\nabla^-+2\ell t_2\nabla^+}Q^{-n})^*$, $\SN_{-t_1,-t_2,n}=\bar Q^{(n)}e^{\frac12rt_1\nabla^--2\ell t_2\nabla^+}$, and ${\SN}_{-t_1,-t_2,n_j}^{\epi(X_0)}$ defined analogously using ${\SN}_{-t_1,-t_2,n_j}$, so that \eqref{eq:bckwdpush1} will follow (after setting $t_1=t_2=t$) if we prove that $(\tfrac{\p}{\p t_1}+\tfrac{\p}{\p t_2})F_{t_1,t_2}(X_0;a)=\cL_{\ell,r}F_{t_1,t_2}(X_0;a)$.
But since $\cL_{\ell,r}=\cL_{\ell,0}+\cL_{0,r}$, it is actually enough to prove that 
\begin{equation}
\tfrac{\p}{\p t_1}F_{t_1,t_2}(X_0;a)=\cL_{0,r}F_{t_1,t_2}(X_0;a)\qand\tfrac{\p}{\p t_2}F_{t_1,t_2}(X_0;a)=\cL_{\ell,0}F_{t_1,t_2}(X_0;a).\label{eq:t1t2}
\end{equation}
The first of these is a version of the Kolmogorov equation for TASEP \eqref{thm:KolmogorovBackwards} where in the contour integrals in \eqref{def:sm} and \eqref{def:sn} there are respectively additional factors $\psi_+(w)$ and $\psi_-(w)$ inside the integrands, both independent of $t$, and it is not hard to see that the proof in Sec. \ref{sec:kolm} works without any difficulty after this modification. We will prove below that \eqref{eq:bckwdpush1} holds in the PushTASEP case $r=0$, and inspecting the proof it is easy again to seen that this implies the second equation in \eqref{eq:t1t2}.
We conclude that it suffices to check the case $r=0$, and we can also take
 $\ell=1$ as it is just a scaling factor.
 
Given this reduction, one may naturally wonder why the proof works for PushASEP and not ASEP.
It just turns out that when one tries to adapt Sch\"utz's formula \cite{MR1468391} to include jumps to the left (for example, as done in \cite{bp-push}) the resulting formula ends up solving the push dynamics instead of the leftward TASEP dynamics.

Having made these comments, we begin the proof.  Proceeding as for TASEP, let $k\in\{1,\dotsc,n\}$ and consider the initial condition $X_0$ and the one where particle $k$ is moved one step to the left along with  its neighboring block.
Write $b$ for the length of this block.
As for TASEP, we need to compute 
\[\big(\SN_{-t,n}^{\epi(\wt X_0)}-\SN_{-t,n}^{\epi(X_0)}\big)(z_1,z_2)
=\ee_{B_0=z_1}\!\left[\SN_{-t,n-\tilde\tau}(B_{\tilde\tau},z_2)\uno{\tilde\tau<n}-\SN_{-t,n-\tau}(B_\tau,z_2)\uno{\tau<n}\right],\]
where now the tildes refer to the initial condition modified in this new way.
We clearly have $\tau=\tilde\tau$ unless $\tilde\tau\in\{k-1,\dotsc,k+b-2\}$ and $B_{\tilde\tau}=X_0(\tau+1)$.
But since the particles $\tilde X_0(j)$ for $j\in\{k,\dotsc,k+b-1\}$ lie in a block, $\tilde\tau\in\{k-1,\dotsc,k+b-2\}$ implies $\tilde\tau=k-1$ (because the walk always jumps down by at least one step).
This means that the above expectation equals
\begin{equation}
\ee_{B_0=z_1}\!\left[\big(\SN_{-t,n-\tilde\tau}(B_{\tilde\tau},z_2)\uno{\tilde\tau<n}-\SN_{-t,n-\tau}(B_\tau,z_2)\uno{\tau<n}\big)\uno{\tilde\tau=k-1,B_{k-1}=X_0(k)}\right]=f_k(z_1)g^{(n)}_{k}(z_2)
\end{equation}
with
\begin{gather}
f_k(z)=\pp_{B_0=z}\!\left(\tilde\tau=k-1,\,B_{k-1}=X_0(k)\right)\\
\shortintertext{and (see also the comment after \eqref{eq:gk} for the meaning of $\tau^{(k)})$}
g^{(n)}_{k}(z)=\Big(\SN_{-t,n-k+1}(X_0(k),z)-\ee_{B_{k-1}=X_0(k)}\!\left[\SN_{-t,n-\tau^{(k)}}(B_{\tau^{(k)}},z)\uno{\tau^{(k)}<n}\right]\Big)\uno{k\leq n}
\end{gather}
(note that $f_k$ does not depend on $b$ but $g^{(n)}_{k}$ does, through $\tau^{(k)}$).
Since this is a rank one kernel, the same arguments as in the TASEP case give
\begin{align}\label{eq:RHSfinal-push}
\cL_{0,1}F_t(X_0;a,n)=-\det(I-\bar\P_aK^\uptext{PushTASEP}_t\bar\P_a)\tr\!\Big[(I-\bar\P_aK^\uptext{PushTASEP}_t\bar\P_a)^{-1}\bar\P_a\big({\textstyle\sum_{k=1}^n}\Delta^{(k)}\big)\bar\P_a\Big]
\end{align}
where $\Delta^{(k)}$ is defined exactly as in \eqref{eq:KcompTASEP} but using the new versions of $f_k$ and $g^{(n)}_k$
(note that here $\cL_{0,1}$ acts on all particles, since all particles are allowed to move in the PushTASEP dynamics).

Now we consider the left hand side. 
As for TASEP, all we need to check in order to finish the proof is that
\begin{equation}\label{eq:toprove-push}
\tfrac{\d}{\d t}K^{\uptext{PushTASEP}}_t={\textstyle\sum_{k\geq1}}\Delta^{(k)}.
\end{equation}
Now we have
\[\Lambda_{i,j}\coloneqq\tfrac{\d}{\d t}K^{\uptext{PushTASEP}}_t(n_i,\cdot;n_j,\cdot)=\big(\tfrac{\d}{\d t}(\SM_{-t,-n_i})^*\big)\SN_{-t,n_j}^{\epi(X_0)}+(\SM_{-t,-n_i})^*(\tfrac{\d}{\d t}\SN_{-t,n_j}^{\epi(X_0)}),\]
and using \eqref{def:sm-push} and \eqref{def:sn-push} we have $\tfrac{\d}{\d t}(\SM_{-t,-n})^*=2(\SM_{-t,-n}\big)^*\nabla^+$ and $\tfrac{\d}{\d t}(\SN_{-t,n})^*=-2\nabla^+(\SM_{-t,n}\big)^*$.
Therefore
\[\Lambda_{i,j}=2\SM_{-t,-n_i}H_{n_j}\]
with
\begin{align}
H_n(z_1,z_2)&=\nabla^+\SN_{-t,n}^{\epi(X_0)}(z_1,z_2)-\ee_{B_0=z_1}\tsm\big[\nabla^+\SN_{-t,n}(B_\tau,z_2)\uno{\tau<n}\big]\\
&=\SN_{-t,n}^{\epi(X_0)}(z_1+1,z_2)-\ee_{B_0=z_1}\tsm\big[\SN_{-t,n-\tau}(B_\tau+1,z_2)\uno{\tau<n}\big].
\end{align}
The first term on the second line can be expressed as $\ee_{B_0=z_1}\tsm\big[\SN_{-t,n-\hat\tau}(B_{\hat\tau}+1,z_2)\uno{\hat\tau<n}\big]$ with $\hat\tau$ the hitting time of the strict epigraph of $\big(X_0(m+1)-1\big)_{m\geq0}$, and thus
\[H_n(z_1,z_2)=\ee_{B_0=z_1}\tsm\big[\SN_{-t,n-\hat\tau}(B_{\hat\tau}+1,z_2)\uno{\hat\tau<n}-\SN_{-t,n-\tau}(B_{\tau}+1,z_2)\uno{\tau<n}\big].\]
Now we have $\hat\tau\leq\tau$, so the difference inside the brackets vanishes for $\hat\tau\geq n$, and it also vanishes when $\hat\tau<n$ and $B_{\hat\tau}\geq X_0(\hat\tau+1)+1$.
Thus $H_n(z_1,z_2)$ equals
\begin{align}
&\sum_{k=0}^{n-1}\pp_{B_0=z_1}\!\left(\hat\tau=k,\,B_{k}=X_0(k-1)\right)\\
&\hspace{0.1in}\times\left(\SN_{-t,n-k}(X_0(k+1)+1,z_2)-\ee_{B_{k}=X_0(k+1)}\!\left[\SN_{-t,n-\tau}(B_{\tau}^{(k)}+1,z_2)\uno{\tau^{(k)}<n}\right]\right)\\
&=\sum_{k=0}^{n-1}\widehat f_{k+1}(z_1)\widehat g^{(n)}_{k+1}(z_2)
\end{align}
with
\[\widehat f_k(z)=\pp_{B_0=z}\!\left(\hat\tau=k-1,\,B_{k-1}=X_0(k)\right)\]
and
\begin{equation}
\widehat g^{(n)}_{k}(z)=\SN_{-t,n-k+1}(X_0(k)+1,z)-\ee_{B_{k-1}=X_0(k)}\!\left[\SN_{-t,n-\tau^{(k)}}(B_{\tau^{(k)}}+1,z)\uno{\tau^{(k)}<n}\right].
\end{equation}
We have then
\[\Lambda_{i,j}=2\sum_{k=1}^{n_j}\SM_{-t,-n_i}\widehat f_k(z_1)\otimes\widehat g^{(n_j)}_{k}(z_2),\]
and so it remains to prove that $\widehat f_k=f_k$ for all $k$ and $\widehat g^{(n)}_{k}=\frac12g^{(n)}_{k}$ for $k\leq n$ (using as for TASEP that $g^{(n)}_k\equiv0$ for $k>n$).

For fixed $k$ let $\tilde\tau_k$ denote the hitting time of the initial condition $\tilde X_0$ introduced above where the block starting at $X_0(k)$ is moved one step to the left, so that $f_k(z)=\pp_{B_0=z}\!\left(\tilde\tau_k=k-1,\,B_{k-1}=X_0(k)\right)$.
We obviously have $\hat\tau\leq\tilde\tau_k$, and thus $\hat\tau=k-1$ and $B_{k-1}=X_0(k)$ implies $\tilde\tau_k=k-1$, so that $\widehat f_k(z)=\pp_{B_0=z}\!\left(\hat\tau=\tilde\tau_k=k-1,\,B_{k-1}=X_0(k)\right)$.
This gives
\[f_k(z)-\widehat f_k(z)=\pp_{B_0=z}\!\left(\tilde\tau_k>\hat\tau=k-1,\,B_{k-1}=X_0(k)\right)=0.\]
Next we prove the identity $\widehat g^{(n)}_{k}=\frac12g^{(n)}_{k}$ for $k\leq n$.
Proceeding as for TASEP, we apply $e^{2t\nabla^+}$ on the right of this identity to see that it is equivalent to
\begin{multline}
\bar Q^{(n-k+1)}(X_0(k)+1,z)-\ee_{B_{k-1}=X_0(k)}\!\left[\bar Q^{(n-\tau^{(k)})}(B_{\tau^{(k)}}+1,z)\uno{\tau^{(k)}<n}\right]\\
=\tfrac12\bar Q^{(n-k+1)}(X_0(k),z)-\tfrac12\ee_{B_{k-1}=X_0(k)}\!\left[\bar Q^{(n-\tau^{(k)})}(B_{\tau^{(k)}},z)\uno{\tau^{(k)}<n}\right].
\end{multline}
Following again the TASEP argument, the two sides are of the form $2^z$ times a polynomial in $z$, so it is enough to prove the equality for $z<X_0(n)$, in which case the identity becomes
\[\pp_{B_{k-1}=X_0(k)}\!\left(B_n=z-1,\,\tau^{(k)}\geq n\right)=\tfrac12\pp_{B_{k-1}=X_0(k)}\!\left(B_n=z,\,\tau^{(k)}\geq n\right),\]
which is easy to prove again since the walk takes Geom$[\frac12]$ steps to the left and $z\leq X_0(n)$.

\vs

\noindent{\bf Acknowledgements.}
MN and JQ were supported by the Natural Sciences and Engineering Research Council of Canada.
DR was supported by Programa Iniciativa Cient\'ifica Milenio grant number NC120062 through Nucleus Millenium Stochastic Models of Complex and Disordered Systems, by Conicyt Basal-CMM Proyecto/Grant PAI AFB-170001, and by Fondecyt Grant 1160174.

\printbibliography[heading=apa]

\end{document}